\documentclass[10pt,twoside]{article}
\usepackage{icmart}%%%%%%%%%%%%%%%%%%
\usepackage{amsmath,amsfonts,amssymb}

\title{\bf Representations of Yangians\vskip 0mm
Associated with Skew Young Diagrams\vskip 6mm}
\author{Maxim Nazarov\hskip-1mm\vspace{-0.5cm}
\thanks{Department of Mathematics, University of York,
York YO10 5DD, England; {\,\tt mln1@york.ac.uk}}}
\date{}
\begin{document}
\maketitle
\thispagestyle{first}\setcounter{page}{1}

%------------------------------------------------------------------------------

\begin{abstract}\vskip2mm

The Yangian of the 
Lie algebra $\mathfrak{gl}_N$ has a distinguished family
of irreducible finite-dimensional representations, called elementary
representations.
They are parametrized by pairs, consisting of a skew Young diagram
and a complex number. Each of these representations has an explicit
realization, it extends the classical realization of the irreducible
polynomial representations of $\mathfrak{gl}_N$ by means of the Young
symmetrizers. We explicitly construct analogues of these elementary
representations for the twisted Yangian, which corresponds to
the Lie algebra $\mathfrak{so}_N$. Our construction provides
solutions to several open problems in the classical representation theory.
In particular, we obtain analogues of the Young symmetrizers for
the Brauer centralizer algebra.

\vskip5mm



\end{abstract}

\vskip8mm

%------------------------------------------------------------------------------

\def\al{\alpha}
\def\AMN{{\rm A}_N(M)}

\def\be{\beta}
\def\BMN{{\rm B}_N(M)}
\def\bx{{\boxed{\phantom{\square}}\kern-.4pt}}

\def\CC{{\mathbb C}}
\def\CNn{(\CC^N)^{\ot\ts n}}
\def\com{\ts,\hskip-.5pt}

\def\de{\delta}
\def\De{\Delta}

\def\End{\operatorname{End}\ts}
\def\Eom{E_{\ts\Om}}

\def\Fom{F_\Om(M)}

\def\g{{\mathfrak{g}}}
\def\ge{\geqslant}
\def\gl{{\mathfrak{gl}}}

\def\id{{\rm id}}

\def\la{\lambda}
\def\La{\Lambda}
\def\lap{\la^{\ts\prime}}
\def\lcd{\ts,\,\ldots,}
\def\le{\leqslant}

\def\mi{{\raise.5pt\hbox{-}}}
\def\mup{\mu^{\ts\prime}}
\def\mv{\kern127pt}
\def\mw{\kern-81pt}

\def\ns{\hskip-1pt}
\def\nup{\nu^{\ts\prime}}

\def\om{\omega}
\def\Om{\Omega}
\def\op{\oplus}
\def\ot{\otimes}

\def\ph{\varphi}

\def\RR{{\mathbb R}} 

\def\sgn{\operatorname{sgn}\ts}
\def\si{\sigma}
\def\sl{\mathfrak{sl}}
\def\so{\mathfrak{so}}
\def\sp{\mathfrak{sp}}
\def\St{\widetilde{S}}

\def\ts{\hskip1pt}
\def\Tt{\ts\widetilde{T\ts}\!\ns}

\def\UM{\operatorname{U}(\gl_M)}
\def\UMN{\operatorname{U}(\gl_{N+M})}
\def\UN{\operatorname{U}(\gl_N)}
\def\US{\operatorname{U}(\so_N)}
\def\USM{\operatorname{U}(\so_M)}
\def\USMN{\operatorname{U}(\so_{N+M})}

\def\Vom{V_{\ts\Om}}
\def\Vt{\widetilde{V}}

\def\Wom{W_\Om(M)}

\def\XMN{\operatorname{X}\ts(\gl_{N+M},\si)}
\def\XN{\operatorname{X}\ts(\gl_N,\si)}

\def\YMN{\operatorname{Y}(\gl_{N+M})}
\def\YN{\operatorname{Y}(\gl_N)}
\def\YS{\operatorname{Y}(\gl_N,\si)}

\def\ZZ{{\mathbb Z}} 

%==============================================================================

\section{Yangian of the general linear Lie algebra}\setzero
\vskip3mm

\noindent\textbf{1.1.}
For each simple finite-dimensional Lie algebra $\g$ over the field $\CC\ts$,
Drinfeld \cite{D1} introduced a canonical deformation of the universal
enveloping algebra of the polynomial current Lie algebra $\g\ts[x]$.
This deformation is a certain Hopf algebra over $\CC\ts$,
denoted by ${\rm Y}(\g)$
and called the \textit{Yangian} of the simple Lie algebra $\g\ts$.
Now consider the general linear Lie algebra $\gl_N\ts$,
it contains the special linear Lie algebra $\sl_N$ as a subalgebra. 
The Hopf algebra which is
called the Yangian of the reductive Lie algebra $\gl_N$ and is
denoted by $\YN$, was considered in the earlier works of
mathematical physicists from St.-Petersburg,
see for instance \cite{FT}.
The Hopf algebra $\YN$ is a deformation of the universal
enveloping algebra of the Lie algebra $\gl_N[x]$,
and the Yangian ${\rm Y}(\sl_N)$ of the simple Lie algebra $\sl_N$
is a Hopf subalgebra of $\YN$. 
Throughout this article, we assume that $N$ is a positive integer.

The unital associative algebra $\YN$ over $\CC$ has a 
family of generators $T_{ij}^{(a)}$ where $a=1,2,\ts\ldots\ts$ and
$i\ts,\ns j=1\lcd N$. The defining relations for these generators
can be written in terms of the formal power series
\begin{equation}\label{1.31}
T_{ij}(x)=
\de_{ij}\cdot1+T_{ij}^{(1)}x^{-\ns1}+T_{ij}^{(2)}x^{-\ns2}+\,\ldots
\,\in\,\YN\,[[x^{-1}]]\,.
\end{equation}
Here $x$ is the formal parameter. Let $y$ be another formal parameter,  
then the defining relations in the associative algebra $\YN$
can be written as
\begin{equation}\label{1.32}
(x-y)\cdot[\ts T_{ij}(x)\ts,T_{kl}(y)\ts]\ts=\;
T_{kj}(x)\ts T_{il}(y)-T_{kj}(y)\ts T_{il}(x)\,,
\end{equation}
where $i\com j\com k\com l=1\lcd N\ts$.
The square brackets in (\ref{1.32}) denote usual commutator.
In terms of the formal series (\ref{1.31}),
the coproduct $\De:\YN\to\YN\ot\YN$ is defined by

\vskip-12pt
\begin{equation}\label{1.33}
\De\bigl(T_{ij}(x)\bigr)\ts=\ts\sum_{k=1}^N\ T_{ik}(x)\ot T_{kj}(x)\,;
\end{equation}

\vskip2pt\noindent
the tensor product on the right hand side of the equality (\ref{1.33})
is taken over the subalgebra $\CC[[x^{-1}]]\subset\YN\,[[x^{-1}]]\ts$.
The counit homomorphism $\varepsilon:\YN\to\CC$ is determined by
the assignment $\,\varepsilon:\,T_{ij}(u)\ts\mapsto\ts\de_{ij}\cdot1$.

For each $i$ and $j$ one can determine a formal power series 
$\Tt_{ij}(x)$ in $x^{-1}$ with the coefficients in $\YN$
and the leading term $\de_{ij}\ts$, by the system of equations
$$
\sum_{k=1}^N\,\ts T_{ik}(x)\ts\Tt_{kj}(x)=\de_{ij}
\quad\textrm{where}\quad
i\com j=1\lcd N.
$$
The antipode ${\rm S}$ on $\YN$ is the anti-automorphism
of the algebra $\YN$, defined by the assignment
$\,{\rm S}:\,T_{ij}(x)\ts\mapsto\ts\Tt_{ij}(x)\ts$.
We also use the involutive automorphism $\xi_N$
of the algebra $\YN$, defined by the assignment
$\,\xi_N:\,T_{ij}(x)\ts\mapsto\ts\Tt_{ij}(-x)\ts$.

Take any formal power series $f(x)\in\CC[[x^{-1}]]$ with the leading
term $1$. The assignment

\vskip-14pt
\begin{equation}\label{1.61}
T_{ij}(x)\mapsto\,f(x)\cdot T_{ij}(x)
\end{equation}

\vskip4pt\noindent
defines an automorphism of the algebra $\YN$,
this follows from (\ref{1.31}) and (\ref{1.32}). 
The Yangian ${\rm Y}(\sl_N)$ is the subalgebra in $\YN$
consisting of all elements, which are invariant under
every automorphism (\ref{1.61}). 

It also follows from  (\ref{1.31}) and (\ref{1.32}) 
that for any $z\in\CC\,$, the assignment
$$
\tau_z:\,T_{ij}(x)\ts\mapsto\,T_{ij}(x-z)
$$
defines an automorphism $\tau_z$ of the algebra $\YN$. Here
the formal power series in $(x-z)^{-1}$ should be re-expanded in $x^{-1}$.
Regard the matrix units $E_{ij}\in\gl_N$ as generators of the universal
enveloping algebra $\UN$. The assignment
$$
\al_N:\ts T_{ij}(x)\,\mapsto\,\de_{ij}\cdot1-E_{ji}\,x^{-1}
$$
defines a homomorphism $\al_N:\YN\to\UN$.
By definition, the homomorphism $\al_N$ is surjective.
For more details and references on the definition of the Yangian $\YN$,
see \cite{MNO}.

\newpage%----------------------------------------------------------------------

\noindent\textbf{1.2.}
Let $\nu=(\nu_{\ts1},\nu_{\ts2}\ts,\,\ldots\,\ts)$ be any partition.
As usual, the parts of $\nu$ are arranged in
the non-increasing order\ts: $\nu_{\ts1}\ge\nu_{\ts2}\ge\ldots\ge0$. Let  
$\nu^{\ts\prime}=(\nup_1,\nup_2\ts,\,\ldots\,\ts)$ be
the partition conjugate to $\nu$. In particular, $\nup_1$
is the number of non-zero parts of the partition $\nu\ts$.
An irreducible module over the Lie algebra $\gl_N$ is called
{\it polynomial\/},
if it is equivalent to a submodule in the tensor product of
$n$ copies of the defining $\gl_N\ts$-module $\CC^N$, 
for some integer $n\ge0$.
The irreducible polynomial $\gl_N\ts$-modules are parametrized by
partitions $\nu$ such that $\nup_1\le N\ts$.
Here $n=\nu_1+\nu_2+\ldots\,\,$.
Let $V_\nu$ be the irreducible module corresponding to $\nu$.
This $\gl_N\ts$-module is of highest weight $(\nu_{\ts1}\lcd\nu_N)$.
Here we choose the  Borel subalgebra in $\gl_N$ consisting of the   
upper triangular matrices, and fix the basis of the diagonal matrix
units $E_{11}\lcd E_{NN}$ in the corresponding Cartan subalgebra of $\gl_N$.

Take any non-negative integer $M$.
Let the indices $i$ and $j$ range over the set $\{1\lcd N+M\}$.
Fix the basis of the matrix units $E_{ij}$ in the Lie algebra
$\gl_{N+M}\ts$. We suppose that the subalgebras $\gl_N$ and $\gl_M$
in $\gl_{N+M}$ are spanned by elements $E_{ij}\ts$ where
respectively $i\com j=1\lcd N$ and $i\com j=N+1\lcd N+M$.
Let $\la$ and $\mu$ be two partitions, such that
$\lap_1\le N+M$ and $\mup_1\le M$. Consider the irreducible
modules $V_\la$ and $V_\mu$ over the Lie algebras $\gl_{N+M}$ and $\gl_M$.
The vector space
\begin{equation}\label{1.0}
{\rm Hom}_{\,\gl_M}(\ts V_\mu\ts\com V_\la\ts)
\end{equation}
comes with a natural action of the Lie algebra $\gl_N\ts$.
This action of $\gl_N$ may be reducible.
The vector space (\ref{1.0}) is non-zero, if and only if $\la_k\ge\mu_k$
and $\lap_k-\mup_k\le N$ for each $k=1\com2\com\ts\ldots$\,; see for
instance \cite{M}\ts.

Denote by $\AMN$ the centralizer of the subalgebra $\UM\subset\UMN$.
The centralizer $\AMN\subset\UMN$ contains $\UN$ as a subalgebra, and
acts naturally in the vector space (\ref{1.0})\ts.
This action is irreducible. 
For every $M$, Olshanski \cite{O1} defined a homomorphism of associative
algebras $\YN\to\AMN$. Along with the centre
of the algebra $\UMN$, the image of this homomorphism generates
the algebra $\AMN$. We use a version of this homomorphism,
it is denoted by $\al_{NM}\ts$.

The subalgebra in $\YMN$
generated by $T_{ij}^{(a)}$ where $i\com j=1\lcd N\,$,
by definition coincides with the Yangian $\YN$.
Denote by $\ph_M$ this natural embedding $\YN\to\YMN$. 
Consider also the involutive automorphism $\xi_{N+M}$
of the algebra $\YMN$.
The image of the homomorphism
$$
\al_{N+M}\circ\ts\xi_{N+M}\circ\ts\ph_M:\ts\YN\to\UMN
$$
belongs to the subalgebra $\AMN\subset\UMN$.
Moreover, this image along with the centre
of the algebra $\UMN$, generates the subalgebra $\AMN$.
For the proofs of these claims, see \cite{MO}.
We use the homomorphism $\YN\to\UMN$
\begin{equation}\label{1.69}
\al_{NM}=\,\al_{N+M}\circ\ts\xi_{N+M}\circ\ts\ph_M\circ\ts\xi_N\ts.
\end{equation}
When $M=0$, the homomorphism (\ref{1.69}) coincides with $\al_N\ts$. 
The intersection of the kernels of all homomorphisms
$\,\al_{\ts N0}\com\al_{\ts N1}\com\al_{\ts N2}\ts,\,\ldots\,$ is zero \cite{O1}.

%------------------------------------------------------------------------------

\vskip\baselineskip\noindent\textbf{1.3.}
The $\AMN$-module (\ref{1.0}) depends on the partitions $\la$ and
$\mu$ via the {\it skew Young diagram}
$$
\om=\{\,(i\com j)\in\ZZ^2\ |\ i\ge1,\ \la_i\ge j>\mu_i\,\}\,.
$$
When $\mu=(0\com0\ts,\ts\ldots\ts)$, this is the usual Young diagram
of the partition $\la$. Consider the $\YN$-module obtained from
the $\AMN$-module (\ref{1.0}) by pulling back through the homomorphism
$\al_{NM}\circ\,\tau_z:\YN\to\AMN$.
Since the central elements of $\UMN$ act in (\ref{1.0}) as scalar operators,
this $\YN$-module is irreducible. It is denoted by $V_\om(z)\ts$, and is
called an \textit{elementary module\/}. Its
equivalence class does not depend on the
choice of the integer $M$, such that $\lap_1\le N+M$ and $\mup_1\le M$.

The elementary modules are distinguished amongst all irreducible
$\YN$-modules by the following theorem. Consider the chain of algebras
\begin{equation}\label{chain}
{\rm Y}(\gl_1)\subset{\rm Y}(\gl_2)\subset\ldots\subset{\rm Y}(\gl_N)\,.
\end{equation}
Here for every $k=1\lcd N-1$ we use the embedding
$\ph_1:{\rm Y}(\gl_k)\to{\rm Y}(\gl_{k+1})\ts$.
Consider the subalgebra of $\YN$ generated by the
centres of all algebras in the chain (\ref{chain}), it is called the
\textit{Gelfand-Zetlin subalgebra\/}. This subalgebra
is maximal commutative in $\YN$; see \cite{C3} and \cite{NO}.
Take any finite-dimensional module $W$
over the Yangian $\YN$.

\vskip.85\baselineskip
{\bf Theorem 1.} 
\textit{Two
conditions on the $\YN$-module\/ $W$ are equivalent\ts:}

{\it a) $W$ is irreducible, and the action of the Gelfand-Zetlin subalgebra
of\/ $\YN$ in\/ $W$ is semi-simple\ts;}

{\it b) $W$ is obtained by pulling back through some automorphism\/ 
{\rm(\ref{1.61})} from the tensor product}

\vskip-16pt
\begin{equation}\label{module}
V_{\om_1}(z_1)\ot\ldots\ot V_{\om_m}(z_m)
\end{equation}

\vskip2pt\noindent
{\it of elementary\/ $\YN$-modules, for some skew Young diagrams\/ 
$\om_1\lcd\om_m$ and for some complex numbers\/ $z_1\lcd z_m$ such that\/
$z_k-z_l\notin\ZZ$ for all\/ $k\neq l$.}

\vskip.85\baselineskip
This characterization of irreducible finite-dimensional
$\YN$-modules with semi-simple action of the Gelfand-Zetlin
subalgebra was conjectured by Cherednik, and was
proved by him \cite{C3} under certain extra conditions
on the module~$W$. In~full generality, Theorem 1 was proved in \cite{NT1}.
An irreducibility criterion for the $\YN$-module {\rm(\ref{module})
with arbitrary parameters $z_1\lcd z_m$ was given in \cite{NT2}.

The classification of all irreducible finite-dimensional $\YN$-modules
has been given by Drinfeld \cite{D2}. However,
the general structure of these modules needs a better
understanding. For instance, the dimensions
of these modules are not explicitly known in general.
The tensor products {\rm(\ref{module})
provide a wide class of irreducible $\YN$-modules,
which can be constructed explicitly.

%------------------------------------------------------------------------------

\vskip\baselineskip\noindent\textbf{1.4.}
The $\YN$-module $V_\om(z)$ has an explicit realization.
It extends the classical
realization of irreducible $\gl_N\ts$-module $V_\nu$ by means
of the Young symmetrizers \cite{W}.

Let us use the standard graphic representation
of Young diagrams on the plane $\RR^2$ with two matrix style coordinates.
The first coordinate increases from top to bottom, the second
coordinate increases from left to right. The element $(i\com j)\in\om$
is represented by the unit box with the bottom right corner
at the point $(i\com j)\in\RR^2$. 

Suppose the set $\om$ consists of $n$ elements.
Consider the \textit{column tableau\/} of shape $\om$.
It is obtained by filling  the boxes of $\om$ with numbers $1\lcd n$
consecutively by columns from left to right, downwards in every column.
Denote this tableau by $\Om$.

\newpage
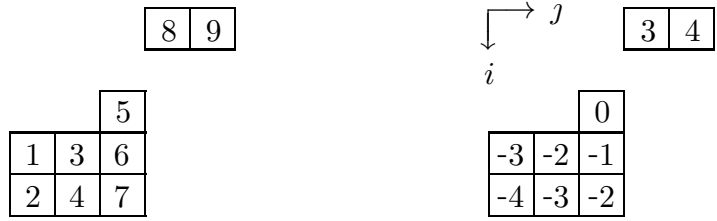
\vbox{
$$
\kern22.4pt\longrightarrow\,j\mw
$$
\vglue-25.2pt
$$
\kern-1pt\vert\mw
$$
\vglue-28pt
$$
\bigr\downarrow\mw\kern0.991pt
$$
\vglue-16pt
$$
i\mw
$$
\vglue-44pt
$$
\phantom{\bx}
\phantom{\bx}
\phantom{\bx}
{\bx}
{\bx}
\kern80pt
\phantom{\bx}
\phantom{\bx}
\phantom{\bx}
{\bx}
{\bx}
$$
\vglue-16.9pt
$$
\phantom{\bx}
\phantom{\bx}
\phantom{\bx}
\phantom{\bx}
\phantom{\bx}
\kern80pt
\phantom{\bx}
\phantom{\bx}
\phantom{\bx}
\phantom{\bx}
\phantom{\bx}
$$
\vglue-16.9pt
$$
\phantom{\bx}
\phantom{\bx}
{\bx}
\phantom{\bx}
\phantom{\bx}
\kern80pt
\phantom{\bx}
\phantom{\bx}
{\bx}
\phantom{\bx}
\phantom{\bx}
$$
\vglue-16.9pt
$$
{\bx}
{\bx}
{\bx}
\phantom{\bx}
\phantom{\bx}
\kern80pt
{\bx}
{\bx}
{\bx}
\phantom{\bx}
\phantom{\bx}
$$
\vglue-16.8pt
$$
{\bx}
{\bx}
{\bx}
\phantom{\bx}
\phantom{\bx}
\kern80pt
{\bx}
{\bx}
{\bx}
\phantom{\bx}
\phantom{\bx}
$$
\vglue-82pt
$$
\kern42pt8\kern9pt9\kern90pt\kern42pt3\kern9pt4\kern-2pt
$$
\vglue-4.5pt
$$
\kern1pt5\kern146pt0
$$
\vglue-16.5pt
$$
1\kern9pt3\kern9pt6\kern116pt\mi3\kern6pt\mi2\kern6pt\mi1\kern26pt
$$
\vglue-16.7pt
$$
2\kern9pt4\kern9pt7\kern116pt\mi4\kern6pt\mi3\kern6pt\mi2\kern26pt
$$
}

\vskip\baselineskip
For each $k=1\lcd n$ put $c_k=j-i$ if the box
$(i,j)\in\om$ is filled with the number $k$ in the tableau $\Om$.
The difference $j-i$ is called the {\it content\/} of the box $(i,j)$
of the diagram $\om$. Our choice of the
tableau $\Om$ provides an ordering of the collection of all contents of $\om$.
In the above figure, on the left we show the column tableau $\Om$
for the partitions $\la=(5,\ns3,\ns3,\ns3,\ns3,\ns0,\ns0,\ts\ldots)$ and
$\mu=(3,\ns3,\ns2,\ns0,\ns0,\ts\ldots)$.
On the right we indicate the contents of all boxes of $\,\om$.

Introduce $n$ complex variables $t_1\lcd t_n$ with the constraints
$t_k=t_l$ for all $k$ and $l$ occuring in the same column of $\Om$\ts.
The number of independent variables among $t_1\lcd t_n$ equals
the number of non-empty columns in the diagram $\om$.
Order lexicographically
the set of all pairs $(k\com l)$ with $1\le k<l\le n$. Take
the ordered product over this set,

\vskip-10pt
\begin{equation}\label{1.3}
\prod_{1\le k<l\le n}^{\longrightarrow}\ 
\left(1-\frac{P_{kl}}{\ts c_k-c_l+t_k-t_l}\ts\right)
\end{equation}

\vskip0pt\noindent
where $P_{kl}$ denotes the operator in the space $\CNn$
exchanging the $k$th and $l$th tensor factors. Consider
(\ref{1.3}) as a function of the constrained variables $t_1\lcd t_n$.

\vskip.85\baselineskip
{\bf Proposition 1.} 
\textit{The rational function\/ {\rm(\ref{1.3})}
is regular at\/ $t_1=\ldots=t_n$.}

\vskip.85\baselineskip
The rational function
(\ref{1.3}) depends only on the differences $t_k-\ts t_l$.
Denote the value of (\ref{1.3}) at $t_1=\ldots=t_n$ by $\Eom\ts$.
Note that for any $\la$ and $\mu\ts$, the linear operator $\Eom$
in the vector space $\CNn$ does not depend on $M$.
For the proof of Proposition~1, see \cite{NT2}.
It provides an explicit expression for the operator~$\Eom\ts$.

Suppose that $\mu=(0\com0\ts,\ts\ldots\ts)\ts$.
In this special case, there is another expression for the operator $\Eom\ts$.
Consider the action of the symmetric group $S_n$ on $\CNn$ by
permutations of the tensor factors. For any $s\in S_n\ts$,
denote by $P_s$ the corresponding operator in $\CNn$. Let $S_\la$
(respectively $S^{\ts\prime}_\la$) be the subgroup in $S_n$ preserving,
as sets, the collections of numbers appearing in every row
(every column) of the tableau $\Om$. Put
$$
X_\la=\sum_{s\in S_\la}\,P_s
\quad\text{and}\quad
Y_\la=\sum_{s\in S^{\ts\prime}_\la}\,P_s\cdot\sgn s
$$
where $\sgn s\ts=\ts\pm1$ is the sign of the permutation $s$.
The product $X_\la Y_\la$ is the \textit{Young symmetrizer\/}
in $\CNn$ corresponding to the tableau $\Om\ts$. We have the equality
\begin{equation}\label{symmetrizer}
\Eom=Y_\la\ts X_\la\ts Y_\la\ts/\,\lap_1!\,\ts\lap_2\ts!\,\ts\ldots\,\ts,
\end{equation}
see \cite{N2}.
In this case, the image of the operator $\Eom$ in $\CNn$ is
equivalent to $V_\la$ as $\gl_N\ts$-module, see \cite{W}.
Here the action of the Lie algebra $\gl_N$ in $\CNn$ is standard.

\newpage%----------------------------------------------------------------------

\noindent\textbf{1.5.}
By pulling the standard action of $\UN$ in the space $\CC^N\ts$
back through the homomorphism

\vskip-16pt
$$
\al_N\circ\,\tau_z:\YN\to\UN\ts,
$$

\vskip3pt\noindent
we obtain a module over the algebra $\YN$, which is denoted by $V(z)$ and
called an \textit{evaluation module}. We have $V(z)=V_\om(z)$
for $\la=(1\com0\ts,\ts\ldots\ts)$ and $\mu=(0\com0\ts,\ts\ldots\ts)$.
For any partitions $\la$ and $\mu$, 
the operator $\Eom$ has the following interpretation, in terms of 
the tensor products of evaluation modules over the Hopf algebra $\YN$.
Let $P_0$ be the operator in $\CNn$ reversing the order of
the tensor factors. 

\vskip.85\baselineskip
{\bf Proposition 2.}
\textit{The operator $\Eom\ts P_0$ is an intertwiner of the $\YN$-modules}
$$
V(c_n+z)\ot\ldots\ot V(c_1+z)
\,\ts\longrightarrow\,
V(c_1+z)\ot\ldots\ot V(c_n+z)\,.
$$

\vskip.15\baselineskip
By Proposition 2, the image of the operator $\Eom$
is a submodule in the tensor product of
evaluation $\YN$-modules $V(c_1+z)\ot\ldots\ot V(c_n+z)$.
Denote this $\YN$-submodule by $\Vom(z)$.
For any $\la$ and $\mu$, we
have the following theorem. Put

\vskip-6pt
\begin{equation}\label{function}
f_\mu(x)\ =\ \prod_{k\ge1}\ 
\frac{(x-\mu_k+k)(x+k-1)}{(x-\mu_k+k-1)(x+k)}\ .
\end{equation}

\vskip-2pt\noindent
This rational function of $x$ expands as a power
series in $x^{-1}$ with the leading~term~$1$.

\vskip.85\baselineskip
{\bf Theorem 2.} 
{\it The\/ $\YN$-module $\Vom(z)\ns$ is equivalent to the elementary
module $V_\om(z)$\/, pulled back through the automorphism of the algebra\/
$\YN$ defined by\/} (\ref{1.61}), {\it where\/} $f(x)=f_\mu(x-z)$.

\vskip.85\baselineskip
This theorem is due to Cherednik \cite{C3}, see also \cite{N}.
It provides an explicit realization of the elementary
$\YN$-module $V_\om(z)$ as a subspace in $\CNn$.
It also shows that the $\YN$-module $\Vom(z)$ is irreducible,
cf.\ \cite{NT2}. The isomorphism between the $\YN$-module $\Vom(z)$,
and the pull-back of the $\YN$-module $V_\om(z)$ as in 
Theorem 1, is unique up to a scalar multiplier.

In Section 2 we give an
analogue of Theorem 2 for the orthogonal Lie algebra $\so_N$,
instead of $\gl_N$.
The case of the symplectic Lie algebra $\sp_N$ is similar to
that of $\so_N$, and is considered in the detailed version
\cite{N} of the present article.

For any simple Lie algebra $\g\ts$ the Yangian ${\rm Y}(\g)$
as defined in \cite{D1}, 
contains the universal enveloping algebra ${\rm U}(\g)$ as a subalgebra.
An embedding $\UN\to\YN$ can be defined by

\vskip-20pt
\begin{equation}\label{embedding}
E_{ij}\ts\mapsto\,-\,T_{ji}^{(1)}.
\end{equation}

\vskip0pt\noindent
The image of ${\rm U}(\sl_N)\subset\UN$ under this emdedding
belongs to ${\rm Y}(\sl_N)\subset\YN$.
The homomorphism $\al_N:\YN\to\UN$
is identical on the subalgebra $\UN$.
The restriction of $\al_N$ to ${\rm Y}(\sl_N)$
provides a homomorphism ${\rm Y}(\sl_N)\to{\rm U}(\sl_N)$,
which is identical on the subalgebra ${\rm U}(\sl_N)$.
For $\g\neq\sl_N$ a homomorphism ${\rm Y}(\g)\to{\rm U}(\g)$
identical on the subalgebra  ${\rm U}(\g)\subset{\rm Y}(\g)$,
does not exist \cite{D1}. For this reason, instead of the Yangian 
${\rm Y}(\so_N)$ from \cite{D1},
we will consider the \textit{twisted Yangian\/}  $\YS$ from \cite{O2}.
Here $\si$ is the involutive automorphism of the Lie algebra $\gl_N$,
such that $-\si$ is the matrix transposition. Then
$\so_N$ is the subalgebra of $\si$-fixed points in $\gl_N\ts$.
%Thus we will actually work with the symmetric pair of Lie algebras 
%$(\ts\gl_N\com\so_N)$.

\newpage%======================================================================

\section{Twisted Yangian of the orthogonal Lie algebra}\setzero
\vskip3mm

\noindent\textbf{2.1.}
The associative algebra $\YS$ is a deformation of the universal
enveloping algebra of the \textit{twisted polynomial current Lie algebra}
$$
\{A(x)\in\gl_N[x]:\sigma(A(x))=A(-x)\}\,.
$$
The deformation $\YS$ is not a Hopf algebra, but a coideal
subalgebra in the Hopf algebra $\YN$.
The definition of the twisted Yangian $\YS$
was motivated by the works of Cherednik \cite{C1} and 
Sklyanin \cite{S} on quantum integrable systems with boundary conditions.
This definition was given by Olshanski in \cite{O2}.

As in Subsection 1.1, let the indices
$i$ and $j$ range over the set $\{1\lcd N\}$. By definition,
$\YS$ is the subalgebra in $\YN$ generated
by the coefficients of all formal power series

\vskip-20pt
\begin{equation}\label{series}
\sum_{k=1}^N\ts\,T_{ki}(-x)\,T_{kj}(x)
\end{equation}
in $x^{-1}$.
Due to (\ref{1.33}),  
the subalgebra $\YS$ in $\YN$ is a right~coideal: 
$$
\De\ts(\ts\YS)\subset\YS\ot\YN\,.
$$

To give the defining relations for the generators of $\YS$, 
introduce the \textit{extended twisted Yangian} $\XN$.
The unital associative algebra $\XN$ has a 
family of generators $S_{ij}^{\ts(a)}$ where $a=1,2,\ts\ldots\,\ts$.
and $i\ts,\ns j=1\lcd N$. Put
\begin{equation}\label{1.771}
S_{ij}(x)=
\de_{ij}\cdot1+S_{ij}^{\ts(1)}x^{-\ns1}+S_{ij}^{\ts(2)}x^{-\ns2}+\,\ldots
\,\in\,\XN\,[[x^{-1}]]\,.
\end{equation}
Defining relations for the generators $\ts S_{ij}^{\ts(a)}\!$
of the algebra $\XN$ can be written~as
$$
(x^2-y^2)\cdot[\ts S_{ij}(x)\ts,S_{kl}(y)\ts]\ts=\;
(x+y)\cdot\bigl(S_{kj}(x)\ts S_{il}(y)-S_{kj}(y)\ts S_{il}(x)\bigr)
$$
$$
-\ (x-y)\cdot\bigl(S_{ik}(x)\ts S_{jl}(y)-S_{ki}(y)\ts S_{lj}(x)\bigr)
+S_{ki}(x)\ts S_{jl}(y)-S_{ki}(y)\ts S_{jl}(x)\,.
$$

\vskip2mm\noindent
All these relations can be written as a single \textit{reflection equation\/},
see \cite{MNO}.
One can define a homomorphism $\pi_N:\XN\to\YS$ by mapping the
series $S_{ij}(x)$ to (\ref{series}).
The homomorpism $\pi_N$ is surjective.
As a two-sided ideal of $\XN$, 
the kernel of the homomorphism $\pi_N$ is generated by
the coefficients of all series
\begin{equation}\label{symmetry}
S_{ij}(x)+(2x-1)\ts S_{ij}(-x)-2x\ts S_{ji}(x)
\end{equation}
in $x^{-1}$.
This ideal is also generated by certain
central elements of $\XN$, see~\cite{MNO}.

The algebra $\XN$ admits
an analogue of the automorphism $\xi_N$ of $\YN$.
Determine a formal power series 
$\St_{ij}(x)$ in $x^{-1}$ with the coefficients in $\XN$
and the leading term $\de_{ij}\ts$, by the system of equations
$$
\sum_{k=1}^N\,\ts S_{ik}(x)\ts\St_{kj}(x)\ts=\,\de_{ij}
\quad\textrm{where}\quad
i\com j=1\lcd N.
$$
Then one can define an involutive automorphism $\eta_N$
of the algebra $\XN$ by the assignment

\vskip-14pt
$$
\textstyle
\eta_N:\ts S_{ij}(x)\mapsto\St_{ij}(-x-\frac{N}2)\,.
$$

\vskip2pt\noindent
However, $\eta_N$ does not determine an automorphism of $\YS$,
because $\eta_N$ does not preserve the ideal of $\XN$ generated by the
coefficients of all series (\ref{symmetry}). 

For any formal power series $f(x)\in\CC[[x^{-1}]]$ with the leading
term $1$, the assignment

\vskip-24pt
\begin{equation}\label{assignment}
S_{ij}(x)\mapsto\,f(x)\cdot S_{ij}(x)
\end{equation}

\vskip-2pt\noindent
defines an automorphism of the algebra $\XN$.
The defining relations of the algebra $\XN$ imply that the assignment

$$
\be_N:\,S_{ij}(x)\,\mapsto\,\de_{ij}\cdot1+\frac{\,E_{ij}-E_{ji}\ts}
{\textstyle x+\frac12}
$$

\vskip0pt\noindent
defines a homomorphism of associative algebras $\be_N:\XN\to\US$.
By definition, the homomorphism $\be_N$ is surjective.
Moreover, $\be_N$ factors through $\pi_N\ts$. Note that the homomorphism
$\YS\to\US$ corresponding to $\be_N\ts$, cannot be obtained from
$\al_N:\YN\to\UN$ by restricting to the subalgebra $\YS$,
because the image of $\YS$ relative to $\al_N$ is not contained in 
the subalgebra $\US\subset\UN$; see \cite{N3}.
An embedding $\US\to\YS$ can be defined by
$$
E_{ij}-E_{ji}\ts\mapsto\,T_{ij}^{\ts(1)}-T_{ji}^{\ts(1)}\,,
$$
cf.\ (\ref{embedding}). 
The homomorphism $\YS\to\US$ corresponding to $\be_N$,
is then identical on the subalgebra $\US\subset\YS$. 

%------------------------------------------------------------------------------

\vskip\baselineskip\noindent\textbf{2.2.}
For any partition $\nu$ with $\nup_1\le N$, the irreducible polynomial 
$\gl_N$-module $V_\nu$ can also be regarded as
a representation of the complex general linear Lie group $GL_N$.
Consider the subgroup $O_N\subset GL_N$ preserving the standard
symmetric bilinear form
$\langle\ ,\,\rangle$ on $\CC^N$. The subalgebra $\so_N\subset\gl_N$
corresponds to this subgroup. Note that the complex Lie group $O_N$
has two connected components.  
In \cite{W} the irreducible finite-dimensional representations of the group
$O_N$ are labeled by the partitions $\nu$ of $n=0\com1\com2\com\,\ldots$
such that  $\nup_1+\nup_2\le N$. Denote by $W_\nu$ 
the  irreducible representation of $O_N$  corresponding to $\nu$.
As $\so_N$-module, $W_\nu$ is irreducible unless
$2\ts\nup_1=N$, in which case $W_\nu$ is a direct sum of two
irreducible $\so_N$-modules.

Choose any embedding of the irreducible representation
$V_\nu$ of the group $GL_N$ into the space $\CNn$. 
Take any two distinct numbers $k,l\in\{1\lcd n\}$.
By applying the bilinear form $\langle\ ,\,\rangle$ to a tensor $w\in\CNn$
in the $k$th and $l$th tensor factors, we obtain a certain tensor
$\widehat{w}\in(\CC^N)^{\ot\ts(n-2)}$. The tensor $w$ is called
\textit{traceless}, if $\widehat{w}=0$ for all distinct $k$ and $l$.
Denote by $\CNn_{\,0}$ the subspace in $\CNn$ consisting of all
traceless tensors, this subspace is $O_N\ts$-invariant.
Then $W_\nu$ can be embedded into $\CNn$ as the intersection
$V_\nu\cap\CNn_{\,0}$, see \cite{W}.

Let the indices $i$ and $j$ range over $\{1\lcd N+M\}$.
Choose the embedding of the Lie algebras $\gl_N$ and $\gl_M$ into
$\gl_{N+M}$ as in Subsection 1.2. It determines embeddings of
groups $GL_N\times GL_M\to GL_{N+M}$ and $O_N\times O_M\to O_{N+M}$.
Take any two partitions $\la$ and $\mu$ such that
$\lap_1+\lap_2\le N+M$ and $\mup_1+\mup_2\le M$.
Consider the irreducible representations $W_\la$ and $W_\mu$ of the groups
$O_{N+M}$ and $O_M$ respectively. The vector space

\vskip-22pt
\begin{equation}\label{2.0}
{\rm Hom}_{\,O_M}(\ts W_\mu\ts\com W_\la\ts)
\end{equation}

\vskip0pt\noindent
comes with a natural action of the group $O_N\ts$.
This action of $O_N$ may be reducible.
The vector space (\ref{2.0}) is non-zero, if and only if $\la_k\ge\mu_k$
and $\lap_k-\mup_k\le N$ for each $k=1\com2\com\ts\ldots$\,; see
\cite{P}\ts. Thus for a given $N$, the vector spaces
(\ref{1.0}) and (\ref{2.0}) are zero or non-zero simultaneously.
Further, for a given
$N$, the dimension of (\ref{2.0}) does not exceed that of (\ref{1.0}).
Our results provide an embedding of (\ref{2.0}) into (\ref{1.0}),
compatible with the action of the orthogonal group $O_N$ in
these two vector spaces.

Denote by $\BMN$ the subalgebra of $O_M$\ts-invariants in the universal
enveloping algebra $\USMN$. Then $\BMN$ contains the subalgebra
$\US\subset\USMN$, and is contained in the centralizer of the subalgebra
$\USM\subset\USMN$. The algebra $\BMN$ naturally acts in the vector
space (\ref{2.0}). The $\BMN$-module (\ref{2.0})
is either irreducible, or splits into a direct
sum of two irreducible $\BMN$-modules. In the latter case,
(\ref{2.0}) is irreducible under the joint action 
of the algebra $\BMN$ and the subgroup $O_N\subset O_{N+M}$.

For every non-negative integer $M$, Olshanski \cite{O2} defined
a homomorphism $\YS\to\BMN$.
Along with the subalgebra of $O_{N+M}\ts$-invariants in
$\USMN$, the image of this homomorphism generates
the algebra $\BMN$. We use a version of this
homomorphism for the algebra $\XN$,
this version is denoted by $\be_{NM}\ts$.

Consider the extended twisted Yangian $\XMN$,
where $-\si$ is the matrix transposition in $\gl_{N+M}$. 
The subalgebra in $\XMN$
generated by $S_{ij}^{(a)}$ where $i\com j=1\lcd N\,$, by definition
coincides with $\XN$.
Denote by $\psi_M$ this natural embedding $\XN\to\XMN$. 
Consider also the involutive automorphism $\eta_{N+M}$
of the algebra $\XMN$. The image of the homomorphism
$$
\be_{N+M}\circ\ts\eta_{N+M}\circ\ts\psi_M:\ts\XN\to\USMN
$$
belongs to the subalgebra $\BMN\subset\USMN$.
Moreover, this image along with the
subalgebra of $O_{N+M}\ts$-invariants in
$\USMN$, generates $\BMN$; see \cite{MO}.
We use the homomorphism $\XN\to\USMN$
\begin{equation}\label{2.69}
\be_{NM}=\,\be_{N+M}\circ\ts\eta_{N+M}\circ\ts\psi_M\circ\ts\eta_N\ts.
\end{equation}
When $M=0$, the homomorphism (\ref{2.69}) coincides with $\be_N\ts$. 
The intersection of the kernels of all homomorphisms
$\,\be_{\ts N0}\com\be_{\ts N1}\com\be_{\ts N2}\ts,\,\ldots\,$
is contained in the kernel of $\pi_N$.

%------------------------------------------------------------------------------

\vskip\baselineskip\noindent\textbf{2.3.}
The $\BMN$-module (\ref{2.0}) depends on the partitions
$\la$ and $\mu$ via the
skew Young diagram $\om$. Using the homomorphism $\be_{NM}:\XN\to\BMN$,
regard (\ref{2.0}) as $\XN$-module. Unlike the
$\YN$-module $V_\om(z)$, this $\XN$-module may depend
on the choice of the integer $M$, such that 
$\lap_1+\lap_2\le N+M$ and $\mup_1+\mup_2\le M$.
Denote this $\XN$-module by $W_\om(M)$.
Note that when $z\neq0$, the automorphism $\tau_z$ of $\YN$ does not
preserve the subalgebra $\YS\subset\YN$. There is no analogue of
the automorphism $\tau_z$ with $z\neq0$ for the algebra $\XN$.

The $O_{N+M}$-invariant elements of $\USMN$ act in (\ref{2.0})
as scalar operators. Thus the $\XN$-module $W_\om(M)$
is either irreducible, or splits into a direct
sum of two irreducible $\XN$-modules. In the latter case, it
becomes irreducible under the joint action 
of the algebra $\XN$ and the subgroup $O_N\subset O_{N+M}$.
Our main result is an explicit realization of the $\XN$-module
$W_\om(M)$, similar to the realization of the elementary
$\YN$-module given by Theorem 2. Our explicit realization is
compatible with the action of the group $O_N$ in $W_\om(M)$.

Take the standard orthonormal basis $e_1\lcd e_N$ in $\CC^N$, so that
$\langle\ts e_i\com e_j\ts\rangle=\de_{ij}\ts$.
The linear operator

\vskip-21pt
\begin{equation}\label{1.45}
u\ot v\,\mapsto\,\langle\ts u\com v\ts\rangle\,\cdot\ts
\sum_{i=1}^N\,\ts e_i\ot e_i
\end{equation}

\vskip-1pt\noindent
in $\CC^N\ns\ot\,\CC^N$
commutes with the action of $O_N$.
Take the complex variables $t_1\lcd t_n$ with the same constraints
as in Proposition 1. Consider the ordered product over the pairs $(k\com l)$,

\vskip-16pt
\begin{equation}\label{1.5}
\prod_{1\le k<l\le n}^{\longrightarrow}\ 
\left(1-\frac{Q_{kl}}
{\ts c_k+c_l+t_k+t_l+N+M}\ts\right)
\end{equation}

\vskip0pt\noindent
where $Q_{kl}$ is the linear operator in $\CNn$,
acting as (\ref{1.45}) in the $k$th and $l$th tensor factors, and acting
as the identity in the remaining $n-2$ tensor factors.
Here the pairs $(k\com l)$ are
ordered lexicographically, as in (\ref{1.3}).
Let us now multiply (\ref{1.5}) by (\ref{1.3}) on the right,
and consider the result as an operator-valued
rational function of the constrained variables $t_1\lcd t_n$.

\vskip.85\baselineskip
{\bf Proposition 3.} 
{\it At\/} $t_1=\ldots=t_n=-\frac12$ {\it the ordered product of\/} 
(\ref{1.5}) {\it and\/} (\ref{1.3}) {\it has the value}

\vskip-10pt
$$
\prod_{(k,\ts l)}^{\longrightarrow}\ 
\left(1-\frac{Q_{kl}}
{\ts c_k+c_l+N+M-1}\ts\right)
\cdot\ts\Eom\,=
$$
\begin{equation}\label{fom}
\Eom\ts\cdot
\prod_{(k,\ts l)}^{\longleftarrow}\ 
\left(1-\frac{Q_{kl}}
{\ts c_k+c_l+N+M-1}\ts\right);
\end{equation}
{\it the ordered products in\/} (\ref{fom}) {\it are taken over all pairs\/}
$(k\com l)$ {\it such that the numbers \/ $k$ and\/ $l$ appear in
different columns of the tableau\/ $\Om$.}

\vskip.85\baselineskip
Denote the operator (\ref{fom}) by $\Fom$.
If $k$ and $l$ appear in different columns of $\Om$,
then

\vskip-16pt
$$
c_k+c_l\ge3-\la'_1-\la'_2\ge3-N-M\,.
$$

\vskip4pt\noindent
Hence each of the denominators in (\ref{fom}) is non-zero for any
choice of $\mu\ts$. The algebra of operators in $\CNn$ generated by all
$P_{kl}$ and $Q_{kl}$ with $1\le k<l\le n\ts$, is called the
\textit{Brauer centralizer algebra\/}; see \cite{B}.
The operator $\Fom$ belongs to this algebra.
Note that the image of the operator $\Fom$ is contained 
in the image of $\Eom\ts$.

Suppose that $M=0$, then $\mu=(0\com0\ts,\ts\ldots\ts)\ts$.
In this special case, the image of the operator $\Eom$ in $\CNn$
is equivalent to $V_\la$ as a representation of the group $GL_N$,
see (\ref{symmetrizer}). It turns out that
the image of the operator $F_\Om(0)$
consists of all traceless tensors from the image of $\Eom\ts$.
In particular, the image of $F_\Om(0)$ in $\CNn$ 
is equivalent to $W_\la$ as a representation of the group $O_N$.
Even in the special case $M=0$, the formulas (\ref{fom}) for
the operator $\Fom$ seem to be new; cf.~\cite{W}.

\newpage%----------------------------------------------------------------------

\noindent\textbf{2.4.}
Let us extend $\si$ to an automorphism of the associative algebra $\UN$.
For any $z\in\CC$, define the \textit{twisted evaluation module\/}
$\Vt(z)$ over the algebra $\YN$
by pulling the standard action of the algebra $\UN$
in the vector space $\CC^N\ts$ back through the 
composition of homomorphisms
$$
\si\ts\circ\,\al_N\circ\,\tau_{-z}\ts:\,\YN\to\,\UN\,.
$$
The evaluation module $V(z)$ and the twisted evaluation module $\Vt(z)$
over $\YN$, have the same restriction to the subalgebra
$\YS\subset\YN\ts$; see (\ref{series}).

For any $\la$ and $\mu$, the operator $\Fom$ has the following
interpretation in terms of the restrictions to $\YS$ of tensor products
of evaluation modules over the Hopf algebra $\YN$; cf.\ Proposition 2.
For each $k=1\lcd n$ put $d_k=c_k+\frac{M}2-\frac12$.
We assume that $\lap_1+\lap_2\le N+M$ and $\mup_1+\mup_2\le M$.

\vskip.85\baselineskip
{\bf Proposition 4.} 
{\it The operator\/ $\Fom$ is an intertwiner of\/ $\YS$-modules}
$$
\Vt(d_1)\ot\ldots\ot\Vt(d_n)
\,\longrightarrow\,
V(d_1)\ot\ldots\ot V(d_n)\,.
$$

\vskip.15\baselineskip
By Proposition 4, the image of $\Fom$
is a submodule in the restriction of the tensor product of
evaluation $\YN$-modules $V(d_1)\ot\ldots\ot V(d_n)$
to  the subalgebra $\YS\subset\YN$.
Denote this $\YS$-submodule by $\Wom$. 
It is also a submodule in the restriction of the $\YN$-module
$V_\Om(\frac{M}2\ns-\ns\frac12)$ to $\YS$. %In the following
%analogue of Theorem 2, we use the rational function (\ref{function}).

\vskip.85\baselineskip
{\bf Theorem 3.} 
{\it a) By pulling the $\XN$-module\/ $W_\om(M)$ back through the 
automorphism of\/ $\XN$ defined by\/} (\ref{assignment}) {\it where\/}
$f(x)=f_\mu(x-\frac{M}2+\frac12)$, 
{\it we get an $\XN$-module that factors through homomorphism\/}
\text{$\pi:\XN\to\YS$.}

{\it b) This\/ $\YS$-module, corresponding to $W_\om(M)$,
is equivalent to\/ $\Wom$.}

\vskip.85\baselineskip
The vector space (\ref{2.0}) of the $\XN$-module 
$W_\om(M)$ comes with a natural
action of the group $O_N$. The action of the group $O_N$ in $\CNn$ 
preserves the image of the operator $\Fom$, because $\Fom$
commutes with this action. Thus the vector space of the $\YS$-module $\Wom$
also comes with an action of $O_N$. The proof of Theorem 3 is
given in \cite{N}. It provides an $O_N$-equivariant isomorphism between
the $\YS$-module corresponding to $W_\om(M)$, and the $\YS$-module $\Wom$.
This isomorphism is unique, up to a scalar multiplier.
The image of the operator $\Fom$ is irreducible
under the joint action of $\YS$ and $O_N$.

Thus we can identify the vector space (\ref{2.0}) with the image
of the operator $\Fom$ uniquely, up to multiplication in (\ref{2.0})
by a non-zero complex number. Using Theorem 2, we can identify
the vector space (\ref{1.0}) with the image of $\Eom\ts$, again uniquely
up to rescaling. Since the image of $\Fom$ is contained in that of
$\Eom\ts$, we then obtain a distinguished embedding of the vector
space (\ref{2.0}) into (\ref{1.0}).

Theorem 3 provides an explicit realization of the $\XN$-module
$W_\om(M)$ as a subspace in $\CNn$. This theorem also
turns the vector space (\ref{2.0}) into a
module over the twisted Yangian $\YS$, equivalent to $\Wom$.
The limited size of the present
article does not allow us to discuss here the analogues of the results
of \cite{NT1} and \cite{NT2} for the $\YS$-modules, obtained in this
particular way; cf.\ \cite{NO}.

\newpage%======================================================================

\end{document}